\documentclass[letterpaper,11pt,oneside,reqno]{Math_class}
\usepackage{float} 
\usepackage[
backend=biber,
style=alphabetic,
sorting=nyt
]{biblatex}
\usepackage{blindtext}
\renewcommand{\dint}{\mathrm{d}}
\title{Inverse of the Gaussian multiplicative chaos: an integration by parts formula}
\author{Tomas Kojar}
\begingroup
\lccode`\~=`\ %
\lowercase{%
  }%
\endgroup
\newcounter{word}

\makeatletter
\newenvironment{proofs}[1][\proofname]{\par
  \pushQED{\qed}%
  \normalfont \topsep0\p@\relax
  \trivlist
  \item[\hskip\labelsep\itshape
  #1\@addpunct{.}]\ignorespaces
}{%
  \popQED\endtrivlist\@endpefalse
}

\hypersetup{
    colorlinks=true,
    linkcolor=blue,
    filecolor=magenta,      
    urlcolor=cyan,
    pdftitle={Inverse of the Gaussian multiplicative chaos : an integration by parts formula},
    pdfpagemode=FullScreen,
    }

\usepackage{fancyhdr}   

\begin{document}
{\let\thefootnote\relax\footnote{{for feedback please contact kojartom@gmail.com}}}

\maketitle
\begin{abstract}
 In this article, we study the analogue of the integration by parts formula from \cite{decreusefond2008hitting} in the context of GMC and its inverse.
\end{abstract}

\tableofcontents

\par\noindent\rule{\textwidth}{0.4pt} 
\setlength{\parindent}{0cm}

\newpage\part{Introduction}
\section{Introduction }
This article is an offshoot application that came up in \cite{binder2023inverse} while doing the preliminary work for extending the work in \cite{AJKS}. In particular, in their work they start with the Gaussian random field $H$ on the circle  with covariance
\begin{equation*}
\Expe{H(z)H(z')}=-\ln\abs{z-z'},    
\end{equation*}
where $z, z'\in\mathbb{C}$ have modulus $1$.  The exponential $\expo{\gamma H}$ gives rise to a random measure $\tau$ on the unit circle $\T$, given by 
\begin{equation*}
\tau(I):=\mu_{H}(I):=\liz{\e}\int_{I}e^{\gamma H_{\e}(x)-\frac{\gamma^{2}}{2}\Expe{\para{H_{\e}(x)}^{2}}}\dx,    
\end{equation*}
for Borel subsets $I\subset \T=\mathbb{R}/ \mathbb{Z}=[0,1)$ and $H_{\e}$ is a suitable regularization. This measure is within the family of \textit{Gaussian multiplicative chaos} measures (GMC) (for expositions see the lectures \cite{robert2010gaussian,rhodes2014gaussian}). So finally, they consider the random homeomorphism $h:[0,1)\to [0,1)$ defined as the normalized measure
\begin{equation*} 
h(x):=\frac{\tau[0,x]}{\tau[0,1]}, x\in [0,1),    
\end{equation*}
and prove that it gives rise to a Beltrami solution and \textit{conformal welding } map. The goal is to extend this result to its inverse $h^{-1}$ and in turn to the composition $h_{1}^{-1}\circ h_{2}$ where $h_{1},h_{2}$ are two independent copies. The motivation for that is of obtaining a parallel point of view of the beautiful work by Sheffield \cite{sheffield2016conformal} of gluing two quantum disks to obtain an SLE loop. \\
We let $Q_{\tau}(x):[0,\tau([0,1])]\to [0,1]$ denote the inverse of the measure $\tau:[0,1]\to [0,\tau([0,1])]$ i.e.
\begin{equation*}
Q_{\tau}(\tau[0,x])=x\tand \tau[0,Q_{\tau}(y)]=y,
\end{equation*}                                  
for $x\in [0,1]$ and $y\in [0,\tau([0,1])].$ The existence of the inverse $Q$ follows from the strict monotonicity of the Liouville measure $\eta$, which in turn follows from being non-atomic \cite[theorem 1]{bacry2003log}. We use the notation Q because the measure $\tau$ can be thought of as the "CDF function" for the "density" $\expo{\gamma H}$ and thus its inverse $\tau^{-1}=Q$ is the quantile (also using the notation $\tau^{-1}$ would make the equations less legible later when we start including powers and truncations). We will also view this inverse as a hitting time for the measure $\tau$
\begin{equation*}
Q_{\tau}(x)=Q_{\tau}(0,x)=T_{x}:=\inf\set{t\geq 0: \tau[0,t]\geq x}.    
\end{equation*}
The inverse homeomorphism map $h^{-1}:[0,1]\to [0,1]$ is defined as
\begin{equation*}\label{inversehomeo}
h^{-1}(x):=Q_{\tau}(x\tau([0,1]))\tfor x\in [0,1]     
\end{equation*}
\noindent Since the inverse of GMC didn't seem to appear in other problems, it was studied very little and so we had to find and build many of its properties. In the article \cite{binder2023inverse}, we go over various basic properties of the inverse $Q$.  Our guide for much for this work was trying to transfer the known properties of the GMC measure to its inverse, the Markovian structure for the hitting times of Brownian motion s (such as the Wald's equation and the independent of the increments of hitting times) and then trying to get whatever property was required for the framework set up by \cite{AJKS} to go through successfully. This was a situation where a good problem became the roadmap for finding many interesting properties for the inverse of GMC and thus GMC itself. \\
When studying the expected value $\Expe{Q(a)}$, we had trouble getting an exact formula. So in the spirit of \cite{decreusefond2008hitting} where they used Malliavin calculus to study the hitting times of processes, we tested using Malliavin calculus to gain better understanding of $\Expe{Q(a)}$. Our guide for applying Malliavin calculus is also the article \cite{aru2020density} where they applied Malliavin calculus to imaginary GMC.
\subsection{Acknowledgements}
We thank I.Binder, Eero Saksman and Antti Kupiainen. We had numerous useful discussions over many years.
\section{Main result}
In \cref{part:IBP}, we study the shifted field $X_{\zeta}=U_{\e}^{r}(\tau_{a}+\zeta)$. We will 
    obtain an integration by parts formula for that field using the techniques from \cite{decreusefond2008hitting,lei2012stochastic,nualart1994skorohod}. Then we will integrate over $\zeta$ to obtain relations for the shifted-GMC and the inverse in \cref{IBPinverse}.
\begin{theorem}
For fixed $\psi\in C_{c}(\Rplus)$ where we normalize $\int_{\R} \psi(a)\da=1$ and $a,L\geq 0$, we have the relation
\begin{eqalign}
\int_{0}^{\infty}\psi(a)\Expe{\eta\para{\tau_{a},\tau_{a}+L}  } \da=&L+\lambda \Expe{ \int_{0}^{r\wedge L}\int_{\zeta}^{\infty}\psi(\eta(\theta-\zeta))\para{\int_{ (\theta-r)\vee 0}^{\theta}\frac{1}{\theta-t}-\frac{1}{r}\deta(t) }   \deta(\theta)\dzeta},
\end{eqalign}
and
\begin{eqalign}
\int_{0}^{\infty}\psi(a)\Expe{\tau_{a}} \da=&\int_{0}^{\infty}\psi(a) a\da+\lambda \Expe{ \int_{0}^{r}\int_{0}^{\infty}\psi(\eta(\theta))\para{\int_{ (\theta+\zeta-r)\vee 0}^{\theta+\zeta}\frac{1}{\theta+\zeta-t}-\frac{1}{r}\deta(t) }   \deta_{\zeta}(\theta)\dzeta},
\end{eqalign}
where $\deta_{\zeta}(\theta):=e^{U(\theta+\zeta)}\dtheta$.
\end{theorem}

\newpage\part{Integration by parts formula}\label{part:IBP}
\section{Setup for Malliavin calculus for the inverse}
In this part we will use the setup from from \cite{decreusefond2008hitting,lei2012stochastic} in order to use the integration by parts formula. In particular, for the Gaussian process $X_{t}:=U_{\epsilon}^{\delta}(t)$ with covariance
\begin{equation}
R(t,s):= \left\{\begin{matrix}
\ln(\frac{r }{\varepsilon} )-\para{\frac{1}{\e}-\frac{1}{r}}\abs{t-s} &, \abs{t-s}\leq \varepsilon\\ 
 \ln(\frac{r}{\abs{t-s}}) +\frac{\abs{t-s}}{r}-1&, \delta>\abs{t-s}\geq \e
\end{matrix}\right.       
\end{equation}
we will use the Malliavin calculus setup for Gaussian processes as developed in \cite{decreusefond2008hitting,lei2012stochastic}. Then once we obtain the various integration by parts formulas, we will then take limit in $\epsilon\to 0$ using the convergence results for GMC (eg.\cite[Theorem 2.1] {berestycki2021gaussian}). For shorthand we will write 
\begin{equation}
\expo{ \bar{U}(t)}=:\expo{\gamma U_{\epsilon}(t)}:=\expo{\gamma U_{\epsilon}(t)-\frac{\gamma^{2}}{2}\ln\frac{1}{\epsilon}}.  \end{equation}
Let $\CH$  be the Hilbert space defined as the closure of the space $\CE$ of step functions on $[0,\infty)$ with respect to the scalar product
\begin{equation}
\braket{1_{[0,s]}  ,1_{[0,t]} }_{\CH}:=R(t,s).    
\end{equation}
The mapping $1_{[0,t]} \mapsto X_{t}$ can be extended to an isometry between $\CH$ and the Gaussian space $H_{1}(X)$ associated with $X$. We will denote this isometry by $\phi\mapsto X(\phi)$. Let $\CS$ be the set of smooth and cylindrical random variables of the form
\begin{equation}
F=f(X(\phi_{1}),...,X(\phi_{n}))    
\end{equation}
for some $n\geq 1$ and $f\in C^{\infty}_{b}(\mathbb{R}^{n})$ (smooth with bounded partial derivatives) and $\phi_{i}\in \CH$. The derivative operator D of a smooth and cylindrical random variable $F\in \CS$ is defined as the $\CH$-valued random variable
\begin{equation}
DF=\sum_{i=1}^{n}\frac{\partial f}{\partial x_{i}} (X(\phi_{1}),...,X(\phi_{n}))  \phi_{i}.   
\end{equation}
The derivative operator $D$ is then a closable operator from $L^2(\Omega)$ into $L^2(\Omega;\CH)$. The Sobolev space $\ud^{1,2}$ is the closure of $\CS$ with respect to the norm
\begin{equation}
 \norm{F}_{1,2}^{2}=E(F^{2})+E(\norm{DF}_{\CH}^{2})   
\end{equation}
The \textit{divergence operator} $\delta$ is the adjoint of the derivative operator. We say that a random variable $u\in L^2(\Omega;\CH)$ belongs to the domain of the divergence operator, denoted by $Dom (\delta)$, if
\begin{equation}
    \abs{E\spara{\braket{DF,u}_{\CH}}}\leq c_{u}\norm{F}_{L^{2}(\Omega)}
\end{equation}
for any $F\in \CS$.  In this case $\delta(u)$ is defined by the duality relationship
\begin{equation}
E\para{F\delta (u)}=  E\spara{\braket{DF,u}_{\CH}},  
\end{equation}
for any $F\in\ud^{1,2}$.

\subsection{Regularity of the covariance}
The following are some of the hypotheses used in the development of Malliavin calculus for Gaussian processes \cite{decreusefond2008hitting,lei2012stochastic}.  The difference is
\begin{equation}
\Expe{\abs{U_{\varepsilon}^{\delta}(t)-U_{\varepsilon}^{\delta}(s)}^{2}}    =2\frac{\abs{t-s}}{\e}(1-\frac{\e}{\delta}),
\end{equation}
which is strictly positive for $t\neq s$. The covariance
\begin{equation}
R(\tau,t):= \left\{\begin{matrix}
\ln(\frac{r }{\varepsilon} )-\para{\frac{1}{\e}-\frac{1}{r}}\abs{\tau-t} &, \abs{\tau-t}\leq \varepsilon\\ 
 \ln(\frac{r}{\abs{\tau-t}}) +\frac{\abs{\tau-t}}{r}-1&, r>\abs{\tau-t}\geq \e
\end{matrix}\right.       
\end{equation}
is in fact an absolutely continuous function as a map $t\mapsto R(\tau,t)$ for each $\tau$: when $\abs{\tau-t}\leq \varepsilon$,  we have the absolutely continuous function $g(t)=\abs{\tau-t}$, and when $\abs{\tau-t}> \varepsilon$, we use that $\ln\frac{1}{x}$ is a differentiable function for $x>\e>0$. We compute the partial derivative to be
\begin{equation}
\ptf{R(\tau,t)}{t}= \left\{\begin{matrix}
-\para{\frac{1}{\e}-\frac{1}{r}}\frac{t-\tau}{\abs{t-\tau}} &, \abs{\tau-t}\leq \varepsilon\\ 
 -\frac{1}{\abs{t-\tau}}\frac{t-\tau}{\abs{t-\tau}} +\frac{1}{r}\frac{t-\tau}{\abs{t-\tau}}&,r>\abs{\tau-t}\geq \e
\end{matrix}\right.       .
\end{equation}
Therefore, for $t>\tau$ the derivative is negative $\ptf{R(\tau,t)}{t}<0$ and for $t<\tau$ it is positive $\ptf{R(\tau,t)}{t}>0$. So it is not continuous on the diagonal, which was one of the constraints in \cite{decreusefond2008hitting}. However, in the work  \cite{lei2012stochastic}, they manage to weaken to the following hypotheses that are satisfied in this setting in \cref{lem:regularitycovariance} 
\begin{lemma}
For all $T>0$ the supremum of the integral of the partial derivative is finite for any $\alpha\geq 1$
\begin{equation}\label{finitint0}
\sup_{s\in [0,T]}\int_{0}^{T}\abs{\ptf{R(s,t)}{t}}^{\alpha}  \dt<\infty  
\end{equation}
and in fact for any continuous function $f$ we have that
\begin{equation}\label{eq:continuousint0}
s\mapsto F(s):=\int_{0}^{T}   f(t)\ptf{R(s,t)}{t}\dt  
\end{equation}
is continuous on $[0,\infty)$.
\end{lemma}
\noindent Finally, because of the stationarity the process $U_{\e}(t)$ does not necessarily diverge to $+\infty$ as $t\to +\infty$. So that means that if we apply the results from \cite{lei2012stochastic}, we have to maintain the upper truncation $\tau_{a}\wedge T$.

\subsection{Regularity of $U_{\epsilon}(\tau_{a})$ and the inverse}
In this section we discuss the Malliavin differentiability situation for $U_{\epsilon}(Q_{\e}(a))$ and for the inverse $Q(x)$, in the limit $\epsilon=0$. For the stopped process there is generally a lack of Malliavin differentiability. For example, for Brownian motion consider any stopping time $T$ eg. the hitting time $T=T_{a}$ of the integrated Geometric Brownian motion of level $a>0$
\begin{equation}
\int_{0}^{T_{a}}e^{B_{s}-\frac{1}{2}s}\ds=a.     
\end{equation}
Then the stopped Brownian motion $W_{T}$ is not Malliavin differentiable (\cite[footnote pg.4]{nam2021locally}). If it was differentiable, we would have that $W_{T}=\int_{0}^{\infty}1_{s\leq T}dW_{s}\in \mathbb{D}^{1,2}$ and $1_{s\leq T}\in \mathbb{D}^{1,2}$. However, by  \cite[proposition 1.2.6]{nualart2006malliavin} we would get that for any $s\geq 0$ either $P[s\leq T]=0$ or $1$, which is a contradiction. \\
On the other hand, for the inverse  for $\epsilon>0$, there are some results. The Malliavin derivative for increasing integral processes has been studied in \cite{nualart1994skorohod}.
\begin{lemma}\cite[lemma 1.5]{nualart1994skorohod}
Let $\set{A_{t}}_{t\in [0,1]}$ be a continuous process such that: 
\begin{enumerate}
    \item Strictly positive $A_{t}>0$ for all $t\in [0,1]$.
    
    \item There exists a version of $A$ such that for all $h\in H$, the map $(\lambda, t)\mapsto A_{t}(\omega+\lambda h)$ is continuous.
    
    \item Finite negative moments $\sup_{t\in [0,1]}A_{t}^{-1}\in L^{p}$ for $p\geq 2$.
    
    \item Finite Malliavin derivative moments: $A\in L^{p}([0,1];\ud^{1,p})$ for $p\geq 2$.
    
\end{enumerate}
For fixed constant $c>0$ consider the hitting time of the integrated process $T_{c}:=\inf\set{t>0: \int_{0}^{t}A_{s}\ds\geq c}$. Then we have $T_{c}\in \ud^{1,p}$ for $p\geq 2$ with Malliavin derivative
\begin{equation}\label{eq:inversederivativeformula}
DT_{c}=\frac{-1}{A_{T_{c}}}\para{\int_{0}^{T}DA_{r}\dr}\ind{T_{c}<1}.     
\end{equation}
\end{lemma}
\noindent In our case we have $A_{t}:=:\expo{\gamma U_{\epsilon}(t)}:$ satisfies all the above assumptions. However, the fraction $\frac{-1}{A_{T_{c}}}=\expo{-\gamma U_{\epsilon}(T_{c})+\frac{\gamma^{2}}{2}\ln\frac{1}{\e}}$ is likely diverging 
because for $c\approx 0$ we have $T_{c}\approx 0$ yet the expectation at zero diverges
\begin{equation}
\Expe{\expo{-\gamma U_{\epsilon}(0)+\frac{\gamma^{2}}{2}\ln\frac{1}{\e}}    }=\expo{\gamma^{2}\ln\frac{1}{\e} }=\e^{-\gamma^{2}}\to +\infty.
\end{equation}
So likely the above formula will not make sense in the limit $\e\to 0$. This lack of differentiability also appears in the works \cite{decreusefond2008hitting,lei2012stochastic}, nevertheless through mollification they manage to extract some interesting formulas that we will try to mimic for the setting of GMC. We apply this first step to the inverse and to match notation write $\tau_{a}:=Q_{\e}(a)$ and also suppress the $\e$ in $\eta(\theta):=\eta_{\e}(\theta)$.\\
We use the same regularization. Suppose that $\phi$ is a nonnegative smooth function with compact support in $(0,+\infty)$ and define for any $T > 0$
\begin{equation}
Y:=\int_{0}^{\infty}\phi(a)\para{\tau_{a}\wedge T }\da.     
\end{equation}
The next result states the differentiability of the random variable $Y$ in the sense of Malliavin calculus and provides an explicit formula for its derivative.
\begin{lemma}\label{derivativemolliinv}
The derivative for the mollified inverse $Y$ is
\begin{align*}
D_{r}Y=&-\gamma\int_{0}^{T}\phi(\eta(\theta))\int_{0}^{\theta}\ind{[0,s]}(r)\deta(s)=-\gamma\int_{\eta(r)}^{\eta(T)}\phi(y)\para{y-\eta(r)}\dtau_{y}. 
\end{align*}
\end{lemma}
\begin{remark}
As we can see in the above formula we get $\dtau_{y}$, which by inverse function theorem is equal to $e^{-\gamma U_{\e}(\tau_{y})+\frac{\gamma^{2}}{2}\ln\frac{1}{\e}}$ in agreement with the formula \cref{eq:inversederivativeformula}.
\end{remark}
\begin{proofs}
Due to $\phi's$ compact support the $Y$ is bounded, and so we can apply Fubini's theorem
\begin{equation}
Y=\int_{0}^{\infty}\phi(a)\int_{0}^{\tau_{a}\wedge T } \dtheta\da    =\int_{0}^{T}\int_{\eta(\theta)}^{\infty}\phi(a)\da \dtheta.
\end{equation}
So here we need to compute the Malliavin derivative of $\eta(\theta)$. By linearity and chain rule for the derivative operator $D$ we obtain
\begin{eqalign}
D_{t}\para{\int_{0}^{x}e^{\gamma U_{\e}(s)-\frac{1}{2}\Exp\spara{\para{\gamma U_{\e}(s)}^{2}}}\ds  }=&    \int_{0}^{x}e^{\gamma U_{\e}(s)-\frac{1}{2}\Exp\spara{\para{\gamma U_{\e}(s)}^{2}}}\gamma D_{t}\para{U_{\e}(s) }\ds\\
=&   \int_{0}^{x}e^{\gamma U_{\e}(s)-\frac{1}{2}\Exp\spara{\para{\gamma U_{\e}(s)}^{2}}}\gamma 1_{[0,s]}(t)\ds\\
=&   \gamma \eta(t,x\vee t ).
\end{eqalign}
Since $\e>0$, we have that $\Expe{\para{\eta(t,x\vee t )}^{2}}<\infty$ and so $\eta(\theta)\in  \ud^{ 1,2}$ (this can also work in the limit $\e=0$ by taking $\frac{2}{\gamma^{2}}>2\Leftrightarrow \gamma<1$). Therefore, by chain rule we get  $Y \in \ud^{ 1,2}$ with
\begin{equation}
D_{r}Y=-\int_{0}^{T}\phi(\eta(\theta))D_{r}(\eta(\theta))\dtheta    =-\int_{0}^{T}\phi(\eta(\theta))\gamma \eta(r,\theta\vee r )\dtheta.
\end{equation}
Finally, making the change of variable $\eta(\theta )= y$ yields
\begin{equation}
D_{r}Y=-\gamma\int_{\eta(r)}^{\eta(T)}\phi(y)\para{y-\eta(r)}\dtau_{y}.
\end{equation}
\end{proofs}

\newpage

\section{Integration by parts formula}

In this section we will obtain an integration by parts formula for $\Expe{\eta(\tau_{a},\tau_{a}+L)}$ using the techniques from \cite{decreusefond2008hitting,lei2012stochastic,nualart1994skorohod}. We apply the Malliavin calculus framework to the Gaussian field $U_{\e_{n}}$ for each fixed $\e_{n}$ and then at the very end we will take limits $\e_{n}\to 0$ in the integration by parts formulas for $\Expe{\eta_{\e_{n}}(\tau_{\e_{n},a},\tau_{\e_{n},a}+L)}$. For simplicity we will temporarily write $\eta=\eta_{\e_{n}}$ and $\tau_{a}=\tau_{\e_{n},a}$.
\subsection{Nonlinear expected value}
\noindent For the usual GMC we know that its expected value is linear $\Expe{\eta(a,b)}=b-a$.  Using the Markovian-like $\delta$-(SMP) property from before, we obtain a nonlinear relation for the expected value of the inverse.
\begin{proposition}\label{differencetermunshifted} We have for $a>0$ and $r\geq \delta $
\begin{eqalign}\label{eq:nonlinearexpect}
\Expe{\eta^{\delta}(Q^{\delta}(a),Q^{\delta}(a)+r)}-r=\Expe{Q^{\delta}(a)}-a&=\int_{0}^{\infty}\Proba{ Q_{R(t)}^{\delta}(a)\leq t \leq  Q^{\delta}(a)}\dt\\
&=\int_{0}^{\infty}\Proba{ \eta^{\delta}(t)\leq a \leq  \eta_{R(t)}^{\delta}(t) }\dt>0.
\end{eqalign}
In particular, for any $a>0$ we have $\Expe{Q^{\delta}(a)}>a$. 
\end{proposition}
\begin{remark}
This proposition shows that the GMC $\eta$ does \textit{not} satisfy a "strong" translation invariance i.e. $\Expe{\eta(Q(a),Q(a)+r)}\neq  r$. So the same is likely true for $Q(a,a+t)$ 
\begin{equation}
\Expe{Q(a,a+t)}=\int_{0}^{\infty}\Proba{t>\eta^{\delta}(Q^{\delta}(a),Q^{\delta}(a)+r)}\dr\neq \int_{0}^{\infty}\Proba{t>\eta^{\delta}(0,r)}=\Expe{Q(t)}.
\end{equation}
It also shows that $\Expe{Q^{\delta}(a)}$ is a nonlinear function of $a$. 
\end{remark}
\begin{remark}
Ideally we would like to check whether the RHS of \cref{eq:nonlinearexpect} is uniformly bounded in $a>0$
\begin{equation}
\supl{a>0}\int_{0}^{\infty}\Proba{\eta(t) \leq a\leq  \eta_{R(t)}(t)}\dt    <\infty\tor =\infty,
\end{equation}
but it is unclear of how the window $[\eta(t) ,\eta_{R(t)}(t)]$ grows as $t\to +\infty$.
\end{remark}
\subsection{Assumptions}
In the work \cite[section 6]{lei2012stochastic,nualart1994skorohod}, they make some assumptions about the covariance $R(s,t)$ of the field that are worth comparing with even though we have to do a new proof for $\eta$.
\begin{itemize}
    \item[(H1)]  For all $t\in [0, T ]$, the map $s\mapsto R(s, t)$ is absolutely continuous on [0, T ] and for some $\alpha>1$ we have
    \begin{equation}
\sup_{s\in [0,T]}\int_{0}^{T}\abs{\ptf{R(s,t)}{t}}^{\alpha}  \dt<\infty.  
\end{equation}

    \item[(H3)] The function $R_t := R(t, t)$ has bounded variation on $[0, T ]$.

    \item[(H5)] $ \limsup_{t\to+\infty} X_t = +\infty$ almost surely.

    \item[(H6)]  For any $0 \leq s < t$, we have 
    \begin{eqalign}
     \Expe{\abs{X_t - X_s}^2 } > 0.   
    \end{eqalign}

    \item[(H7)]  For any continuous function $f$ , we have that
\begin{equation}
s\mapsto F(s):=\int_{0}^{T}   f(t)\ptf{R(s,t)}{t}\dt  
\end{equation}
is continuous on $[0,\infty)$.
    
\end{itemize}
Even though our setting is different since we study hitting times of $\eta(t)$ and not of $X_{t}$, these assumptions have analogues. In the  \cref{sec:regcovariance} we compute the derivative of
\begin{equation}
R(\tau,t):= \left\{\begin{matrix}
\ln(\frac{r }{\varepsilon} )-\para{\frac{1}{\e}-\frac{1}{r}}\abs{\tau-t} &, \abs{\tau-t}\leq \varepsilon\\ 
 \ln(\frac{r}{\abs{\tau-t}}) +\frac{\abs{\tau-t}}{r}-1&, r>\abs{\tau-t}\geq \e
\end{matrix}\right.       
\end{equation}
to be 
\begin{equation}
\ptf{R(\tau,t)}{t}= \left\{\begin{matrix}
-\para{\frac{1}{\e}-\frac{1}{r}}\frac{t-\tau}{\abs{t-\tau}} &, \abs{\tau-t}\leq \varepsilon\\ 
 -\frac{1}{\abs{t-\tau}}\frac{t-\tau}{\abs{t-\tau}} +\frac{1}{r}\frac{t-\tau}{\abs{t-\tau}}&,r>\abs{\tau-t}\geq \e
\end{matrix}\right.       .
\end{equation}
and show the assumptions $(H1),(H3)$ and $H(7)$. The assumption $(H6)$ is immediate from the covariance computation. Finally the analogue of the assumption $(H6)$ for $\eta$ is immediate since it is in fact a strictly increasing function.
\subsection{Integration by parts formula for truncated hitting time}
As in these works here too we study the exponential evaluated at the stopping time: 
\begin{equation}
M_{t+\zeta}:=\expo{\lambda U_{\e_{n}}^{\delta}(t+\zeta)- \frac{\lambda^{2}}{2}\ln\frac{1}{\e_{n}}   },
\end{equation}
$\tfor t,\zeta\geq 0$ and some $\lambda\in [0,\sqrt{2})$. The $\zeta$ is important here because we will then integrate over $\zeta$ to obtain a formula  for $\Expe{\eta(\tau_{a},\tau_{a}+L)}$ with $a,L\geq 0$.
\noindent The following proposition follows from \cite[prop.2.1]{decreusefond2008hitting} and it asserts that $\delta_{t}M:=\frac{1}{\lambda}\para{M_{t+\zeta}-1}$ satisfies an integration by parts formula, and in this sense, it coincides with an extension of the Skorokhod divergence of $M \one_{[0,t]}$. 
\begin{proposition}\label{dualityrelation}\cite[prop.2.1]{decreusefond2008hitting}
For any smooth and cylindrical random variable of the form $F=f(X_{t_{1}},...,X_{t_{n}})$ for $t_{i}\in [0,t]$, we have
\begin{equation}\label{Lipintpartsdisc}
\Expe{F\delta_{t}M}=\Expe{\sum_{i=1}^{n}\frac{\partial f}{\partial x_{i}}(X_{t_{1}},...,X_{t_{n}})\int_{0}^{t+\zeta}M_{s}\ptf{R}{s}(s,t_{i})\ds   }.
\end{equation}
\end{proposition}
\noindent By writing
\begin{equation}\label{eq:regularizedinverse}
Y=\int_{0}^{\infty}\phi(a)\para{\tau_{a}\wedge T}\da=\int_{0}^{T}\int_{\eta(\theta)}^{\infty}\phi(a)\da\dtheta,
\end{equation}
where $\phi\in C_{c}^{\infty}(\Rplus)$, we will apply  \cref{dualityrelation} to $F:=p(Y-t)$, where $p\in C_{c}^{\infty}(\R)$ and $M_{t+\zeta}$. In particular, due to the discontinuity of the $\ptf{R}{s}$ along the diagonal, we choose $p_{\delta}(x-y)=0$ when $x>y$ as they do in \cite[Theorem 6.5.]{lei2012stochastic}. The following lemma uses the proof structure of \cite[lemma 6.4]{lei2012stochastic}.
\begin{lemma}
We have the integration by parts relation
\begin{eqalign}\label{lem:MYintparts}
\Expe{p(Y)\delta_{t}M}=&- \Expe{p'(Y)\int_{0}^{T}\phi(\eta(\theta) )  \int_{0}^{t+\zeta}M_{s}\spara{\int_{0}^{\theta}    \ptf{R}{s}(b,s)\deta(b)  }  \ds \dtheta} .   \\
=&- \Expe{p'(Y)\int_{0}^{\eta(T)}\phi(y)  \int_{0}^{t+\zeta}M_{s}\spara{\int_{0}^{y}    \ptf{R}{s}(\tau_{b},s)\db  }  \ds \dtau_{y}} .
\end{eqalign}
\end{lemma}
\begin{remark}
The inverse $\tau_{y}$ is a strictly increasing continuous function (even at the limit $\e=0$) and so we can define its Riemann-Stieltjes integral. This is because of the a)non-atomic nature of GMC \cite[theorem 1]{bacry2003log} and b)GMC;s continuity and strict monotonicity, which in turn follows from satisfying bi-\Holder over dyadic intervals \cite[theom 3.7]{AJKS}.
\end{remark}
\begin{proofs}
The strategy is to discretize the domain $[0,T]$ and thus bring us to the setting of proposition \cref{dualityrelation}. Consider an increasing sequence $D_N:=\set{\sigma_{i}: 0=:\sigma_{0}<\sigma_{1}<...<\sigma_{N}:=T}$ of finite subsets of $[0,T]$ such that their union $\bigcup_{N\geq 1}D_{N}$ is dense in $[0,T]$. Set $D_{N}^{\theta}:=D_N\cap [0,\theta]$ with $\sigma(\theta):=\max(D_{N}^{\theta})$, to let
\begin{equation}
\eta_{N}(\theta):=\eta_{N}(\sigma(\theta)):=\sum_{k=1}^{\sigma(\theta)}\expo{\bar{U}_{\e}(\sigma_{k})}    \para{\sigma_{k}-\sigma_{k-1}}
\end{equation}
and 
\begin{equation}
Y_N:= \int_{0}^{T}\psi(\eta_{N}(\theta) )  \dtheta=\sum_{m=1}^{N}\psi(\eta_{N}(\sigma_{k-1}) ) \para{\sigma_{k}-\sigma_{k-1}}.   
\end{equation}
 Then, $Y_N$ and $p(Y_N)$ are Lipschitz functions of $\set{U_{\e}(t) :t \in D_N}$. The partial $\sigma_{i}$-derivative is
 \begin{equation}
\frac{\partial (p(Y_N))}{\partial \sigma_{i}}=-p'(Y_{N})\sum_{k=i+1}^{N}\phi(\eta_{N}(\sigma_{k-1}) )\cdot\para{ \expo{\bar{U}_{\e}(\sigma_{i})}    \para{\sigma_{i}-\sigma_{i-1}} }\cdot\para{\sigma_{k}-\sigma_{k-1}}         
 \end{equation}
 and so the formula \cref{dualityrelation} implies that
 \begin{eqalign}
\Expe{p(Y_{N})\delta_{t}M}=&- \Expe{\sum_{i=2}^{N}p'(Y_{N})\sum_{k=i+1}^{N}\phi(\eta_{N}(\sigma_{k-1}) )\cdot\para{ \expo{\bar{U}_{\e}(\sigma_{i})}    \para{\sigma_{i}-\sigma_{i-1}}} \cdot\para{\sigma_{k}-\sigma_{k-1}} \para{\int_{0}^{t+\zeta}M_{s}\ptf{R}{s}(\sigma_{i},s) \ds } } \\
=&- \Expe{p'(Y_{N})\sum_{k=2}^{N}\phi(\eta_{N}(\sigma_{k-1}) )  \int_{0}^{t+\zeta}M_{s}\spara{\sum_{i=1}^{k-1}\expo{\bar{U}_{\e}(\sigma_{i})}    \ptf{R}{s}(\sigma_{i},s)\para{\sigma_{i}-\sigma_{i-1}}}  \ds \para{\sigma_{k}-\sigma_{k-1}} }.
 \end{eqalign} 
The function $r \mapsto \int_{0}^{t+\zeta}M_{s}\ptf{R}{s}(s,r) \ds$ is continuous and bounded by condition (H1). As a consequence, we can take the $N$-limit of the above Riemann sum to get the integral formula
\begin{equation}
\Expe{p(Y)\delta_{t}M}=- \Expe{p'(Y)\int_{0}^{T}\phi(\eta(\theta) )  \int_{0}^{t+\zeta}M_{s}\spara{\int_{0}^{\theta} \expo{\bar{U}_{\e}(b)}    \ptf{R}{s}(b,s)\db  }  \ds \dtheta} .    
\end{equation}
Finally, making the change of variable $\eta(\theta )= y$ yields
\begin{eqalign}
\Expe{p(Y)\delta_{t}M}=&- \Expe{p'(Y)\int_{0}^{\eta(T)}\phi(y)  \int_{0}^{t+\zeta}M_{s}\spara{\int_{0}^{\tau_{y}} \expo{\bar{U}_{\e}(b)}    \ptf{R}{s}(b,s)\db  }  \ds \dtau_{y}}    \\
=&- \Expe{p'(Y)\int_{0}^{\eta(T)}\phi(y)  \int_{0}^{t+\zeta}M_{s}\spara{\int_{0}^{y}      \ptf{R}{s}(\tau_{b},s)\db  }  \ds \dtau_{y}}    ,
\end{eqalign}
where in the last equality we used that $\eta$ and $\tau$ are inverses of each other.
\end{proofs}
\subsection{Limits in the Integration by parts relation}
\noindent In this section we set a specific regularization $\phi_{\e}(x)=\frac{1}{\e}\one_{[-1,0]}(\frac{x}{\e})$ in \cref{eq:regularizedinverse}
\begin{equation}
    Y_{\e,a}:=\int_{0}^{\infty}\phi_{\e}(x-a)(\tau_{x}\wedge T)\dx=\frac{1}{\e}\int_{a-\e}^{a}(\tau_{x}\wedge T)\dx=\int_{0}^{1}(\tau_{a-\e\xi}\wedge T)\dxi,
\end{equation}
where we let $\tau_{x}=0$ when $x<0$, and we take limits of $\phi=\phi_{\e}$ and $p=p_{\delta}$ as $\e,\delta\to 0$. Before that step, since the derivative of the mollification $p'$ will diverge in the limit $\delta\to 0$, we first integrate both sides in \cref{lem:MYintparts} as done in \cite[Proof of Proposition 6.1]{lei2012stochastic}.
\begin{proposition}\label{prop:IBPformulaUtaua}
Fix $\psi\in C_{c}^{\infty}(\Rplus)$  and set $c:=\int_{\R} \psi(a)\da$. We have the following integration by parts relation
\begin{eqalign}\label{intbypartsM}
&\int_{0}^{\infty}\psi(a)\int_{0}^{\infty}\Expe{p_{\delta}(Y_{\e,a}-t)M_{t+\zeta}  } \dt\da\\
=&c-\lambda   \Expe{ \int_{0}^{\infty}\int_{0}^{\eta(T)}\para{\int_{0}^{1}\psi(y+\e w)p_{\delta}(Y_{\e,y-\e w}-t)\dw}  M_{t+\zeta}\spara{\int_{0}^{y}    \ptf{R}{t}(\tau_{b},t+\zeta)\db  }}.
\end{eqalign}
By further taking the limits in $\e,\delta\to 0$ we obtain the following relation for each $T>0$
\begin{eqalign}\label{intbypartsM2}
\int_{0}^{\infty}\psi(a)\Expe{M_{\tau_{a}\wedge T+\zeta}  } \da=c-\lambda   \Expe{ \int_{0}^{\eta(T)}\psi(y)M_{\tau_{y}+\zeta}\spara{\int_{0}^{y}    \ptf{R}{t}(\tau_{b},\tau_{y}+\zeta)\db  }   \dtau_{y}}.
\end{eqalign}
\end{proposition}
\begin{remark}\label{rem:shiftedGMCcorr}
By integrating over $\zeta\in [0,L]$ we obtain an IBP for shifted-GMC
\begin{eqalign}\
&\int_{0}^{\infty}\psi(a)\Expe{\eta\para{\tau_{a}\wedge T,\tau_{a}\wedge T+L}  } \da\\
=&c (L-0)-\lambda   \Expe{ \int_{0}^{L}\int_{0}^{\eta(T)}\psi(y)M_{\tau_{y}+\zeta}\spara{\int_{0}^{y}    \ptf{R}{t}(\tau_{b},\tau_{y}+\zeta)\db  }   \dtau_{y}\dzeta}.
\end{eqalign}
\end{remark}

\begin{proofs}
Continuing from \cref{lem:MYintparts} we rewrite it as 
\begin{eqalign}
\int_{0}^{\infty}\Expe{p_{\delta}(Y_{\e,a}-t)M_{t+\zeta}  } \dt=&1+\lambda    \int_{0}^{\infty}\Expe{p_{\delta}(Y_{\e,a}-t)\delta(M\one_{[0,t+\zeta]}  }\dt    \\
=&1-\lambda    \int_{0}^{\infty}\Expe{p_{\delta}'(Y_{\e,a}-t)\int_{0}^{\eta(T)}\phi_{\e}(y-a) \int_{0}^{t+\zeta}M_{s}\spara{\int_{0}^{y}    \ptf{R}{s}(\tau_{b},s)\db  }  \ds \dtau_{y} }\dt.
\end{eqalign}
Now to remove the $p'$ issue, we do an integration by parts for the $\dt$ integral to obtain
\begin{eqalign}
1-\lambda    \int_{0}^{\infty}\Expe{p_{\delta}(Y_{\e,a}-t)\int_{0}^{\eta(T)}\phi_{\e}(y-a)  M_{t+\zeta}\spara{\int_{0}^{y}    \ptf{R}{t}(\tau_{b},t+\zeta)\db  }   \dtau_{y} }\dt.
\end{eqalign}
We multiply both sides by $\psi(a)$ and integrate over the variable $a$
\begin{eqalign}
&\int_{\R}\psi(a)\int_{0}^{\infty}\Expe{p_{\delta}(Y_{\e,a}-t)M_{t+\zeta}  } \dt\da\\
=&c-\lambda   \Expe{\int_{\R}\psi(a) \para{\int_{0}^{\infty}p_{\delta}(Y_{\e,a}-t)\int_{0}^{\eta(T)}\phi_{\e}(y-a)  M_{t+\zeta}\spara{\int_{0}^{y}    \ptf{R}{t}(\tau_{b},t+\zeta)\db  }   \dtau_{y}\dt}\da}.
\end{eqalign}
Here for the $\da$-integral we use that $\phi_{\e}(y-a)=\frac{1}{\e}\one_{[-1,0]}(\frac{y-a}{\e})$ to write
\begin{eqalign}
&c-\lambda   \Expe{ \int_{0}^{\infty}\int_{0}^{\eta(T)}\para{\frac{1}{\e}\int_{y}^{y+\e}\psi(a)p_{\delta}(Y_{\e,a}-t)\da}  M_{t+\zeta}\spara{\int_{0}^{y}    \ptf{R}{t}(\tau_{b},t+\zeta)\db  }   \dtau_{y}\dt}.
\end{eqalign}
Finally, we do a change of variable $a = y +\e w$
\begin{eqalign}\label{eq:maintermwithlimitstoswap}
&c-\lambda   \Expe{ \int_{0}^{\infty}\int_{0}^{\eta(T)}\para{\int_{0}^{1}\psi(y+\e w)p_{\delta}(Y_{\e,y-\e w}-t)\dw}  M_{t+\zeta}\spara{\int_{0}^{y}    \ptf{R}{t}(\tau_{b},t+\zeta)\db  }   \dtau_{y}\dt}\\
=:&c-\lambda   \Expe{ \int_{0}^{\infty}\int_{0}^{\eta(T)}F_{\e,\delta}(y,t)   G(t,y)\dtau_{y}\dt}
\end{eqalign}
for 
\begin{eqalign}
F_{\e,\delta}(y,t) :=&\int_{0}^{1}\psi(y+\e w)p_{\delta}(Y_{\e,y-\e w}-t)\dw,\\
 G(t,y):=&M_{t+\zeta}\spara{\int_{0}^{y}    \ptf{R}{t}(\tau_{b},t+\zeta)\db  }.    
\end{eqalign}
We next take limits and justify their swapping with the integrals.
\proofparagraph{Limit $\e\to 0$}
 We use that the inverse $\tau_{y}$ is a continuous function to take limit
\begin{equation}
\liz{\e}Y_{\e,y-\e w}=\liz{\e}\int_{0}^{1}(\tau_{y-\e w-\e\xi}\wedge T)\dxi  = \tau_{y}\wedge T
\end{equation}
and so the limiting $w$-integral is
\begin{eqalign}
&\liz{\e}\int_{0}^{1}\psi(y+\e w)p_{\delta}(Y_{\e,y-\e w}-t)\dw\\
=&\int_{0}^{1}\psi(y)p_{\delta}(\tau_{y}w+\tau_{y} (1-w)-t)\dw\\
=&\psi(y)p_{\delta}(\tau_{y}-t).
\end{eqalign}
We next justify that we can swap limit and integrals in \cref{eq:maintermwithlimitstoswap}. By the compact support and smoothness of $\phi\tand p$ we have a uniform constant
\begin{equation}
F_{\e,\delta}(y,t)=\int_{0}^{1}\psi(y+\e w)p_{\delta}(Y_{\e,y-\e w}-t)\dw\leq K.
\end{equation}
Moreover,we can assume that compact support is contained $\text{supp}p_{\delta}\subseteq [0,T+\delta]$ and so the infinite integral in \cref{lem:regularitycovariance} gets restricted to $[0,T+\delta]$. We also use the uniform constant to bound as follows
\begin{eqalign}\label{lem:regularitycovariancetwo}
 \eqref{lem:regularitycovariance} \leq  K\Expe{ \int_{0}^{T+\delta}\int_{0}^{\eta(T)}   \abs{G(t,y)}\dtau_{y}\dt}. 
\end{eqalign}
Finally, we will need to revert to the previous formula in terms of GMC
\begin{eqalign}
\int_{0}^{y}      \ptf{R}{s}(\tau_{b},s)\db=\int_{0}^{\tau_{y}}    \ptf{R}{s}(b,s)\deta(b).
\end{eqalign}
We put all these together
\begin{eqalign}
  &\Expe{ \int_{0}^{\infty}\int_{0}^{\eta(T)}F_{\e,\delta}(y,t)   G(t,y)\dtau_{y}\dt}\\
  \leq &K\Expe{ \int_{0}^{\eta(T)}   \dtau_{y}\int_{0}^{T+\delta}\para{\int_{0}^{T+\delta}   \abs{ \ptf{R}{t}(b,t+\zeta)}\deta(b)} M_{t+\zeta}\dt}\\
  = &KT\Expe{\int_{\zeta}^{\zeta+T+\delta}\int_{0}^{T+\delta}   \abs{ \ptf{R}{t}(b,t)}\deta(b)\deta(t)}\\
   = &KT\int_{\zeta}^{\zeta+T+\delta}\int_{0}^{T+\delta}   \abs{ \ptf{R}{t}(b,t)}\db\dt,
\end{eqalign}
where we also used that $\tau_{y}\leq T+\delta$ and applied Fubini-Tonelli to integrate-out the GMCs. This final quantity is indeed finite due to the continuity of the integral as explained in \cref{lem:regularitycovariance}. Therefore, all together we can use dominated convergence theorem to swap limits and integral
\begin{eqalign}\label{eq:maintermwithlimitstoswapfordelta}
\liz{\e}\eqref{eq:maintermwithlimitstoswap}=&c-\lambda   \Expe{ \int_{0}^{T+\delta}\int_{0}^{\eta(T)}\psi(y)p_{\delta}(\tau_{y}-t)  M_{t+\zeta}\spara{\int_{0}^{y}    \ptf{R}{t}(\tau_{b},t+\zeta)\db  }   \dtau_{y}\dt}.
\end{eqalign}
\proofparagraph{Limit $\delta\to 0$}
Here we follow parts of the \cite[proof of lemma 3.3]{decreusefond2008hitting}.  Here we just use from \cref{lem:regularitycovariance} that the integral
\begin{equation}
\int_{0}^{y}    \ptf{R}{t}(\tau_{b},t+\zeta)\db=\int_{0}^{\tau_{y}} \expo{\bar{U}_{\e}^{r}(b)}    \ptf{R}{t}(b,t+\zeta)\db    
\end{equation}
is continuous in $t$ even if $\zeta=0$ but as long as $\e_{n}>0$. Therefore, we can take the limit in $\delta\to 0$. Now in terms of using dominated convergence theorem, we use the same dominating factor as above.\\
In summary we get the following limit
\begin{eqalign}
\liz{\delta}\eqref{eq:maintermwithlimitstoswapfordelta}=&c-\lambda   \Expe{ \int_{0}^{\eta(T)}\psi(y)M_{\tau_{y}+\zeta}\spara{\int_{0}^{y}    \ptf{R}{t}(\tau_{b},\tau_{y}+\zeta)\db  }   \dtau_{y}}.
\end{eqalign}

\end{proofs}
\section{Formula for the shifted GMC}
In this section we use the IBP formula in \cref{prop:IBPformulaUtaua} to obtain a formula for the shifted GMC and the expected value of the hitting time. We will work with field $U_{\varepsilon}^{r}$ for $r>\e>0$ and $\zeta>0$. As mentioned in \cref{rem:shiftedGMCcorr} we already have one formula. By integrating over $\zeta\in [0,L]$ we obtain an IBP for shifted-GMC
\begin{eqalign}\label{eq:IBPformulawithpsiandepsilon}
&\int_{0}^{\infty}\psi(a)\Expe{\eta\para{\tau_{a}\wedge T,\tau_{a}\wedge T+L}  } \da\\
=&c L-\lambda   \Expe{ \int_{0}^{L}\int_{0}^{\eta(T)}\psi(y)M_{\tau_{y}+\zeta}\spara{\int_{0}^{y}    \ptf{R}{t}(\tau_{b},\tau_{y}+\zeta)\db  }   \dtau_{y}\dzeta}.
\end{eqalign}
In the rest of the section we try to simplify this formula. 
\subsection{Limit in $\e\to 0$ for fixed $\psi$}
In the \cref{eq:IBPformulawithpsiandepsilon}, ideally one would like to investigate  taking $\e\to 0$ and having the support of the $\psi=\psi_{n}$ to be approximating to a point $a_{0}$. Assuming one can swap limits with integrals one would get the following formula
\begin{eqalign}
&\Expe{\eta\para{\tau_{a_{0}}\wedge T,\tau_{a_{0}}\wedge T+L}  } \da\\
=&c L-\lambda   \Expe{ \int_{0}^{L}M_{\tau_{a_{0}}+\zeta}\spara{\int_{0}^{a_{0}}    \ptf{R}{t}(\tau_{b},\tau_{a_{0}}+\zeta)\db  }   \frac{1}{M_{\tau_{a_{0}}}}\dzeta},
\end{eqalign}
where the factor $\frac{1}{M_{\tau_{a_{0}}}}$ originated from the \textit{formal} limit of $\frac{\dtau_{y}}{\dy}=e^{-\bar{U}_{\e}^{r}(\tau_{y})}$. The issue here is that this latter limit doesn't exist because the normalization is reversed (the same is true even for the field $e^{-\bar{U}_{\e}^{r}(s)}$ over deterministic $s$ since its mean is diverging like $\e^{-\gamma^{2}}$.)\\
Therefore, we will study the IBP formula for fixed $\psi$ and $\e\to 0$.
\begin{proposition}\label{prop:fixedpsilimitepsilon}
For fixed $\psi\in C_{c}(\Rplus)$ where we normalize $\int_{\R} \psi(a)\da=1$, we have the relation
\begin{eqalign}\label{eq:ibpformulaforeacht}
\int_{0}^{\infty}\psi(a)\Expe{\eta\para{\tau_{a}\wedge T,\tau_{a}\wedge T+L}  } \da=&L+\lambda \Expe{ \int_{0}^{r}\int_{\zeta}^{\zeta+T}\psi(\eta(\theta-\zeta))\spara{\int_{ (\theta-r)\vee 0}^{\theta}\frac{1}{\theta-t}-\frac{1}{r}\deta(t) }   \deta(\theta)\dzeta},
\end{eqalign}
where the GMCs have the field with $\e=0$.  For simplicity we take $T\geq 1>\e>0$.
\end{proposition}
\begin{remark}
One corollary is the inequality
\begin{eqalign}
\int_{0}^{\infty}\psi(a)\Expe{\eta\para{\tau_{a}\wedge T,\tau_{a}\wedge T+L}  } \da\geq L.    
\end{eqalign}
Here we can actually take limit of $\psi=\psi_{n}$ whose support is converging to a fixed value $a_{0}$, to get the inequality
\begin{eqalign}
\Expe{\eta\para{\tau_{a_{0}},\tau_{a_{0}}+L}  } \geq L,   \end{eqalign}
which agrees with the result in \cref{differencetermunshifted}.
\end{remark}

\subsubsection{Proof of \cref{prop:fixedpsilimitepsilon}}
\noindent  We start by writing the IBP formula explicitly using the covariance function.
\begin{lemma}\label{lemma:covarianceepsilon}
Using the explicit formula of the covariance we have the expression
\begin{eqalign}
&\int_{0}^{y}    \ptf{R}{t}(\tau_{b},\tau_{y}+\zeta)\db  =-\int_{a}^{b}\frac{1}{\tau_{y}+\zeta-t}-\frac{1}{r}\deta(t) -\para{\frac{1}{\e}-\frac{1}{r}} \etam{b, \tau_{y}},
\end{eqalign}
for
\begin{eqalign}
a:=  \tau_{y}\wedge (\tau_{y}+\zeta-r)\vee 0\tand b:=  \tau_{y}\wedge (\tau_{y}+\zeta-\e)\vee 0.    
\end{eqalign}
\end{lemma}
\begin{proofs}
For ease of notation in the proof we let $s:=\tau_{y}+\zeta$ and
\begin{equation}
a:=  \tau_{y}\wedge (s-r)\vee 0, b:=  \tau_{y}\wedge (s-\e)\vee 0, c:= \tau_{y}\wedge (s+\e)\tand d:= \tau_{y}\wedge (s+r).
\end{equation}
Using the explicit formula for the partial derivative in \cref{eq:covariancederivative} we have the following
\begin{eqalign}
&\int_{0}^{\tau_{y}} \expo{\bar{U}_{\e}^{r}(t)}    \ptf{R}{s}(t,s)\dt\\
=&\int_{a}^{b}\frac{-1}{s-t}\deta(t)+ \frac{1}{r}  \etam{a,b} + \int_{c}^{d}\frac{1}{t-s}\deta(t)+ \frac{-1}{r}   \etam{c,d}+\para{\frac{1}{\e}-\frac{1}{r}} \para{\etam{s\wedge \tau_{y},c}-\etam{b,s\wedge \tau_{y}}}.
 \end{eqalign}
For $s=\tau_{y}+\zeta$ we have
 \begin{equation}
a=  \tau_{y}\wedge (\tau_{y}+\zeta-r)\vee 0, b=  \tau_{y}\wedge (\tau_{y}+\zeta-\e)\vee 0, c:= \tau_{y}\tand d:= \tau_{y}.
\end{equation}
Therefore, the above simplifies
\begin{eqalign}
&\int_{0}^{\tau_{y}} \expo{\bar{U}_{\e}(t)}    \ptf{R}{s}(t,s)\mid_{s=\tau_{y}+\zeta}\dt\\
=&\int_{a}^{b}\frac{-1}{s-t}\deta(t)+ \frac{1}{r}  \etam{a,b} +0+ \frac{-1}{r}   \cdot 0+\para{\frac{1}{\e}-\frac{1}{r}} \para{0-\etam{b, \tau_{y}}}\\
=&-\int_{a}^{b}\frac{1}{s-t}-\frac{1}{r}\deta(t) -\para{\frac{1}{\e}-\frac{1}{r}} \etam{b, \tau_{y}}.
 \end{eqalign}
\end{proofs}
\noindent Returning to \cref{eq:IBPformulawithpsiandepsilon} we write
\begin{eqalign}\label{eq:IBPformulawithpsiandepsilonexplicit}
\eqref{eq:IBPformulawithpsiandepsilon}=&\int_{0}^{\infty}\psi(a)\Expe{\eta\para{\tau_{a}\wedge T,\tau_{a}\wedge T+L}  } \da\\
=&c L-\lambda   \Expe{ \int_{0}^{L}\int_{0}^{\eta(T)}\psi(y)M_{\tau_{y}+\zeta}\spara{-\int_{a}^{b}\frac{1}{\tau_{y}+\zeta-t}-\frac{1}{r}\deta(t) -\para{\frac{1}{\e}-\frac{1}{r}} \etam{b, \tau_{y}}  }   \dtau_{y}\dzeta}\\
=&c L-\lambda   \Expe{ \int_{0}^{L}\int_{0}^{T}\psi(\eta(\theta))M_{\theta+\zeta}\spara{-\int_{\wt{a}}^{\wt{b}}\frac{1}{\theta+\zeta-t}-\frac{1}{r}\deta(t) -\para{\frac{1}{\e}-\frac{1}{r}} \etam{\wt{b}, \theta}  }   \dtheta\dzeta},
\end{eqalign}
where we also undid the change of variables $\tau_{y}=\theta\doncl y=\eta(\theta)$, and let
\begin{eqalign}
\wt{a}:=  \theta\wedge (\theta+\zeta-r)\vee 0\tand \wt{b}:=  \theta\wedge (\theta+\zeta-\e)\vee 0.    
\end{eqalign}
\noindent Taking $\e\to 0$ on the LHS is clear since $\psi$ is compactly supported and bounded. The question is what happens in the RHS. We study each term.
\begin{lemma}\label{lem:zerolimit}
We have the limit
\begin{eqalign}
\liz{\e}\Expe{ \int_{0}^{L}\int_{0}^{T}\psi(\eta(\theta))M_{\theta+\zeta}\spara{ -\para{\frac{1}{\e}-\frac{1}{r}} \etam{\wt{b}, \theta}  }   \dtheta\dzeta}=0.    
\end{eqalign}
\end{lemma}
\begin{proofs}
In the term $\etam{\wt{b}, \theta}$, since $\wt{b}:=  \theta\wedge (\theta+\zeta-\e)\vee 0 $, we have that as soon as $\zeta\geq \e$, we get identically zero $\etam{\wt{b}, \theta}   =0$ for every $\e> 0$. So we just study the integrals
\begin{eqalign}
&\Expe{ \int_{0}^{\e}\int_{0}^{T}\psi(\eta(\theta))M_{\theta+\zeta}\spara{ -\para{\frac{1}{\e}-\frac{1}{r}} \etam{(\theta+\zeta-\e)\vee 0, \theta}  }   \dtheta\dzeta}\\
=&-\para{\frac{1}{\e}-\frac{1}{r}} \int_{0}^{\e}\Expe{ \int_{\zeta}^{\zeta+T}\psi(\eta(\theta-\zeta))\etam{(\theta-\e)\vee 0, \theta-\zeta}    \deta(\theta)}\dzeta.
\end{eqalign}
Here we can apply Lebesgue differentiation theorem. We study the difference of functions
\begin{eqalign}
f(\zeta)-g_{\e}(\zeta):=   \Expe{ \int_{\zeta}^{\zeta+T}\psi(\eta(\theta-\zeta))\etam{0, \theta-\zeta}    \deta(\theta)}-\Expe{   \int_{\e}^{\zeta+T}\psi(\eta(\theta-\zeta))\etam{0,\theta-\e}  \deta(\theta)}.
\end{eqalign}
In the first function by taking limit $\e\to 0$ we get
\begin{eqalign}\label{eq:firstfunctionlimitldt}
\dashint_{0}^{\e}   f(\zeta)\dzeta\to f(0)=   \Expe{ \int_{0}^{T}\psi(\eta(\theta))\etam{ 0,\theta}    \deta(\theta)}.
\end{eqalign}
In the second function, we separate the two limits
\begin{eqalign}\label{eq:secondfunctionlimitldt}
&\dashint_{0}^{\e}\Expe{   \int_{\zeta}^{\zeta+T}\psi(\eta(\theta-\zeta))\etam{0,\theta}  \deta(\theta)}\dzeta+\dashint_{0}^{\e}\Expe{   \int_{\e}^{\zeta+T}\psi(\eta(\theta-\zeta))\para{\etam{0,\theta-\e}-\etam{0,\theta}}  \deta(\theta)}\dzeta.
\end{eqalign}
The first term converges to the same limit as in \cref{eq:firstfunctionlimitldt} and so they cancel out. Therefore, it suffices to show that the second term in \cref{eq:secondfunctionlimitldt} goes to zero. We pull out the supremum 
\begin{eqalign}
&\abs{\dashint_{0}^{\e}\Expe{   \int_{\e}^{\zeta+T}\psi(\eta(\theta-\zeta))\para{\etam{0,\theta-\e}-\etam{0,\theta}}  \deta(\theta)}\dzeta}\\
\leq &\dashint_{0}^{\e}\Expe{\sup_{\e\leq z\leq \e+T}\abs{\etam{z-\e,z}} \cdot  \int_{\zeta}^{\zeta+T}\psi(\eta(\theta-\zeta))  \deta(\theta)}\dzeta.
\end{eqalign}
The quantity inside the expectation is uniformly bounded in $\e$ because we can use \Holder to separate them
\begin{eqalign}
 \Expe{\sup_{\e\leq z\leq \e+T}\abs{\etam{z-\e,z}}^{2}}^{1/2}\cdot\Expe{\para{ \int_{\zeta}^{\zeta+T}\psi(\eta(\theta-\zeta))  \deta(\theta)}^{2}}^{1/2},   
\end{eqalign}
where due to \cref{prop:maxmoduluseta} the first factor goes to zero as $\e\to 0$.
\end{proofs}
\noindent We return to take the limit $\e\to 0$ in \cref{eq:IBPformulawithpsiandepsilonexplicit} 
\begin{eqalign}\label{eq:IBPformulawithpsiandepsilonexplicitfirstterm}
\liz{\e}\eqref{eq:IBPformulawithpsiandepsilonexplicit}=&\int_{0}^{\infty}\psi(a)\Expe{\eta\para{\tau_{a}\wedge T,\tau_{a}\wedge T+L}  } \da\\
=&c L-\lambda  \liz{\e} \Expe{ \int_{0}^{L}\int_{\zeta}^{\zeta+T}\psi(\eta(\theta-\zeta))\para{-\int_{\wt{a}}^{\wt{b}}\frac{1}{\theta-t}-\frac{1}{r}\deta(t) }   \deta(\theta)\dzeta},
\end{eqalign}
for
\begin{eqalign}
\wt{a}:=  \para{\theta-\zeta}\wedge (\theta-r)\vee 0\tand \wt{b}:=  \para{\theta-\zeta}\wedge (\theta-\e)\vee 0.    
\end{eqalign}
We note here that if $\zeta\geq r$, then we get $\wt{a}=\theta-\zeta=\wt{b}$ and so the inner integral becomes zero. So we are left with
\begin{eqalign}
&\liz{\e} \Expe{ \int_{0}^{r}\int_{\zeta}^{\zeta+T}\psi(\eta(\theta-\zeta))\para{-\int_{ (\theta-r)\vee 0}^{ \wt{b}}\frac{1}{\theta-t}-\frac{1}{r}\deta(t) }   \deta(\theta)\dzeta}.
\end{eqalign}
\noindent The following lemma concludes the proof of \cref{prop:fixedpsilimitepsilon}.
\begin{lemma}\label{lemma:limitzeroendpoint}
We have the limit
\begin{eqalign}
&\liz{\e} \Expe{ \int_{0}^{r}\int_{\zeta}^{\zeta+T}\psi(\eta(\theta-\zeta))\para{\int_{ (\theta-r)\vee 0}^{\para{\theta-\zeta}\wedge (\theta-\e)\vee 0}\frac{1}{\theta-t}-\frac{1}{r}\deta(t) }   \deta(\theta)\dzeta}\\
=&\Expe{ \int_{0}^{r}\int_{\zeta}^{\zeta+T}\psi(\eta(\theta-\zeta))\para{\int_{ (\theta-r)\vee 0}^{\theta}\frac{1}{\theta-t}-\frac{1}{r}\deta(t) }   \deta(\theta)\dzeta}.
\end{eqalign}
\end{lemma}
\begin{remark}
A a heuristic we study the integrals without any GMCs:
\begin{eqalign}
\int_{0}^{r}\int_{\zeta}^{\zeta+T}\para{\int_{ (\theta-r)\vee 0}^{\theta-\zeta}\frac{1}{\theta-t}-\frac{1}{r}\dt }   \dtheta\dzeta=&\int_{0}^{r}\int_{0}^{T}\para{\int_{ (\theta+\zeta-r)\vee 0}^{\theta}\frac{1}{\theta+\zeta-t}\dt }   \dtheta\dzeta-r\para{T-\frac{r}{6}}\\
=&\int_{0}^{r}\int_{0}^{T} \ln\frac{1}{\zeta}-\ln\frac{1}{r\wedge \para{\theta+\zeta}}   \dtheta\dzeta-r\para{T-\frac{r}{6}}\\
=&-r\ln\frac{1}{r}\para{1-\frac{3r}{2}}-r\para{T-\frac{r}{6}}.
\end{eqalign}
So we see that even for $r\to 0$ we still have finiteness in the limit $\e\to 0$.  \end{remark}
\begin{proofs}[proof of \cref{lemma:limitzeroendpoint}]
We will apply dominated convergence theorem. In terms of limits we study the inner integrals
\begin{eqalign}
f(\zeta):=\Expe{\int_{\zeta}^{\zeta+T}\psi(\eta(\theta-\zeta))\para{\int_{ (\theta-r)\vee 0}^{\para{\theta-\zeta}\wedge (\theta-\e)\vee 0}\frac{1}{\theta-t}-\frac{1}{r}\deta(t) }   \deta(\theta)}   
\end{eqalign}
Since we have fixed $\psi$ and it has compact support, we get that it is bounded and so we upper bound
\begin{eqalign}
f(\zeta)\leq &K \Expe{\int_{\zeta}^{\zeta+T}\para{\int_{ (\theta-r)\vee 0}^{\para{\theta-\zeta}\wedge (\theta-\e)\vee 0}\frac{1}{\theta-t}\deta(t) }   \deta(\theta)}\\
=&K\int_{\zeta}^{\zeta+T} \int_{ (\theta-r)\vee 0}^{\para{\theta-\zeta}\wedge (\theta-\e)\vee 0}\frac{1}{\para{\theta-t}^{1+\gamma^{2}}}    \dtheta\\
\lessapprox & \frac{T}{\zeta^{\gamma^{2}}},  
\end{eqalign}
where we evaluate the correlation for the two GMCs. This factor is still integrable as long as $\gamma^{2}<1$. Therefore, we can indeed apply the dominated convergence theorem.    
\end{proofs}

\newpage\subsubsection{IBP Formula for inverse}\label{IBPinverse}
We justify taking infinite limit $T\to +\infty$.
\begin{lemma}\label{lem:ibpinfinitet}
The finite $T$ limit of \cref{eq:ibpformulaforeacht} is
 \begin{eqalign}
\int_{0}^{\infty}\psi(a)\Expe{\eta\para{\tau_{a},\tau_{a}+L}  } \da=&L+\lambda \Expe{ \int_{0}^{r}\int_{0}^{\infty}\psi(\eta(\theta))\spara{\int_{ (\theta+\zeta-r)\vee 0}^{\theta+\zeta}\frac{1}{\theta+\zeta-t}-\frac{1}{r}\deta(t) }   \deta_{\zeta}(\theta)\dzeta},
\end{eqalign}
where we used the notation $\deta_{\zeta}(\theta):=e^{U^{r}_{0}(\zeta+\theta)}\dtheta$.
\end{lemma}
\noindent Therefore, for $L\geq r$ we use to \cref{differencetermunshifted}  obtain the following formula for the expected value of the inverse.
\begin{corollary}
The inverse satisfies the following integration by parts formula
\begin{eqalign}
\int_{0}^{\infty}\psi(a)\Expe{\tau_{a}} \da=&\int_{0}^{\infty}\psi(a) a\da+\lambda \Expe{ \int_{0}^{r}\int_{0}^{\infty}\psi(\eta(\theta))\spara{\int_{ (\theta+\zeta-r)\vee 0}^{\theta+\zeta}\frac{1}{\theta+\zeta-t}-\frac{1}{r}\deta(t) }   \deta_{\zeta}(\theta)\dzeta}.
\end{eqalign}
\end{corollary}
\begin{proofs}[proof of \cref{lem:ibpinfinitet}]
Since $\psi$ is compactly supported $supp(\psi)\subset [0,S]$ for some $S>0$ we get that the integral is zero as soon as
\begin{eqalign}
\eta(\theta)>S.    
\end{eqalign}
So for the LHS in \cref{eq:ibpformulaforeacht} we have
\begin{eqalign}
\int_{0}^{S}\psi(a)\Expe{\eta\para{\tau_{a}\wedge T,\tau_{a}\wedge T+L}  } \da.    
\end{eqalign}
Since the shifted GMC $\eta\para{\tau_{a}\wedge T,\tau_{a}\wedge T+L}$ is continuous and uniformly bounded in $T$
\begin{eqalign}
\eta\para{\tau_{a}\wedge T,\tau_{a}\wedge T+L}\leq \eta\para{0,\tau_{a}+L},    
\end{eqalign}
we can apply dominated convergence theorem. For the RHS we start by undoing the change of variables $\theta\leftrightarrow \tau_{y}$ to write
\begin{eqalign}
\eqref{eq:ibpformulaforeacht}=L+\lambda \Expe{ \int_{0}^{r}\int_{0}^{\eta(T)\wedge S}\psi(y)\spara{\int_{ (\tau_{y}+\zeta-r)\vee 0}^{\tau_{y}+\zeta}\frac{1}{\tau_{y}+\zeta-t}-\frac{1}{r}\deta(t) }   e^{U_{\tau_{y}+\zeta}}\dy\dzeta}.    
\end{eqalign}
Here we use the following limiting ergodic statements for GMC \cite[lemma 1]{allez2013lognormal}.
\begin{lemma}\label{ergodiclemma}
Let $M$ be a stationary random measure on $\R$ admitting a moment of order $1+\delta$ for $\delta>0$. There is a nonnegative integrable random variable $Y\in L^{1+\delta}$ such that, for every bounded interval $I\subset \R$, $$\lim_{T \to \infty} \frac{1}{T} M\left(T I\right) = Y |I|\quad \text{almost surely and in }L^{1+\delta},$$
where  $|\cdot|$ stands for the Lebesgue measure on $\R$. As a consequence, almost surely the random measure $$A\in \mathcal{B}(\R)\mapsto \frac{1}{T}M(TA)$$ weakly converges towards $Y|\cdot|$ and $\E_Y[M(A)]=Y |A|$ ($\E_Y[\cdot]$ denotes the conditional expectation with respect to $Y$).
\end{lemma}
\noindent For GMC  the $Y$ variable is equal to one $Y=1$. One way to see it is using the independence of distant GMCs. By splitting $\frac{\eta^{1}(0,n)}{n}$ into alternating even and odd intervals $[k,k+1]$ to get two independent sequences and then apply strong law of large numbers to get convergence to $\frac{\eta^{1}(0,n)}{n}\stackrel{a.s.}{\to}\frac{1}{2}\Expe{\eta^{1}(0,1)}+\frac{1}{2}\Expe{\eta^{1}(1,2)}=1$.\\
Therefore, since the quantity is uniformly bounded in $T$ by bounding by the integral over $\int_{0}^{S}$, we can apply dominated convergence theorem. 
   
\end{proofs}

\newpage\part{Further directions and Appendix}
\section{Further research directions }\label{furtherresearchdirections}
\begin{enumerate}
    
   \item\textbf{ Joint law for the Liouville measure}\\
   The density of the inverse is in terms of two-point joint law of GMC:
   \begin{equation*}
     \Proba{b\geq Q(x)\geq a}=\Proba{\eta(b)\geq x\geq \eta(a)}.  
   \end{equation*}
   (Of course, if we have differentiability, we can just study  $\dfrac{\dint}{\db}\Proba{\eta(b)\geq x}$).  The same issue showed up when studying the decomposition of the inverse.  For example, we could turn the conditional moments' bounds into joint law statements by rewriting the event $\set{Q(a)-Q(b)=\ell}$ in terms of $\eta$. Some approaches include conformal field theory in \cite{zhu2020higher} and \textit{possibly} Malliavin calculus \cite{airault2010stokes,takeuchi2016joint}. See here for work on GMC and Malliavin calculus \cite{broker2020geometry}. It would also be interesting to get bounds on the single and joint density of GMC using the Malliavin calculus techniques in \cite{nourdin2009density}. In the same spirit as in \cite{wong2020universal}, one can try to Goldie-renewal result: see \cite{damek2021diagonal} for recent work extending the Goldie renewal result used in \cite{wong2020universal} to the case of joint law.

    \item \textbf{Regularity for GMC's Malliavin derivative }
It would be interesting to explore the regularity of the Malliavin derivative $D^{k}\eta$ for $k=k(\gamma)$ as $\gamma\to 0$. This can give different upper bounds for the density:
\begin{proposition}
\label{densityestimates}
Let $q, \alpha, \beta$ be three positive real numbers such that $\frac{1}{q}+\frac{1}{\alpha}+\frac{1}{\beta}=1$. Let $F$ be a random variable in the space $\ud^{2,\alpha}$, such that $\Expe{\norm{DF}_{H}^{-2\beta}} < \infty$. Then the density $p(x)$ of $F$ can be estimated as follows
\begin{equation*}
p(x)\leq c_{q, \alpha, \beta}\para{\Proba{\abs{F}>\abs{x}}}^{1/q}  \spara{\Expe{\norm{DF}_{H}^{-1}}+\norm{D^{2}F}_{L^{\alpha}(\Omega;H\otimes H)}  \Expe{\norm{DF}_{H}^{-2\beta}}^{1/\beta}  },  
\end{equation*}
where $\norm{u}_{L^{\alpha}(\Omega;H\otimes H)}:= \Expe{\norm{u}_{H\otimes H}^{\alpha}}^{1/\alpha}$.
\end{proposition}

\item\textbf{Derivatives in the IBP-formula} In the spirit of the derivative computations done in  \cite{decreusefond2008hitting}, one could try to extract some pdes/odes. We included some some heuristics computations for
\begin{eqalign}
M_{\tau_{a_{0}} }:=e^{\lambda U_{\tau_{a_{0}}}-\frac{\lambda^{2}}{2}\ln\frac{1}{\e}}.    
\end{eqalign}
In \cref{intbypartsM2}, we can concentrate $\psi$ around the point $a_{0}$ and use \cref{lemma:covarianceepsilon} to get the identity
\begin{eqalign}
\Psi(a,\lambda):=\Expe{M_{\tau_{a_{0}} }  }    =& 1+\lambda   \Expe{ M_{\tau_{a_{0}} } \int_{ (\tau_{a_{0}} -r)\vee 0}^{ (\tau_{a_{0}} -\e)\vee 0}\para{\frac{1}{\tau_{a_{0}} -t}-\frac{1}{r}}\deta(t)  }\\
&+\lambda \para{\frac{1}{\e}-\frac{1}{r}}  \Expe{ M_{\tau_{a_{0}} }  \etam{ (\tau_{a_{0}} -\e)\vee 0,\tau_{a_{0}}  }}.
\end{eqalign}
So the $\lambda$ derivative of the LHS is: 
\begin{eqalign}
\ptf{\Psi(a,\lambda)}{\lambda}=&\Expe{M_{\tau_{a} }U_{\e}(\tau_{a})  }-\frac{\lambda}{2}\Expe{M_{\tau_{a} } }\ln\frac{1}{\e}\\
=&\frac{1}{\lambda}\Expe{M_{\tau_{a} }\ln\para{M_{\tau_{a}}}}
\end{eqalign}
and of the RHS is
\begin{eqalign}
\Expe{\ptf{\psi(a,\lambda)}{\lambda}}=\ptf{\Psi(a,\lambda)}{\lambda}=&\Expe{M_{\tau_{a}}F(a)}+\lambda \Expe{M_{\tau_{a} }U_{\e}(\tau_{a})F(a)}\\
&-\frac{\lambda^{2}}{2}\Expe{M_{\tau_{a} }F(a) }\ln\frac{1}{\e}\\
=&\Expe{M_{\tau_{a}}\para{1+\ln\para{M_{\tau_{a}}}} F(a)}\\
=&\Expe{\psi(\lambda,a) F(a)}+\Expe{\psi(\lambda,a)\ln\para{\psi(\lambda,a)} F(a)},
\end{eqalign}
where $\psi(\lambda,a):=M_{\tau_{a} }$. So one ODE from here is
\begin{equation}
y'=y(1+\ln(y))c , \twith y(0)=1  
\end{equation}
which has the unique solution 
\begin{equation}
y(\lambda)=\expo{\expo{\frac{c^{2}}{2}}-1}     .
\end{equation}
The $\Psi$ itself satisfies
\begin{eqalign}
\ptf{\Psi(a,\lambda)}{\lambda}=&\Expe{M_{\tau_{a}}F(a)}+\lambda \Expe{M_{\tau_{a} }U_{\e}(\tau_{a})F(a)}\\
&-\frac{\lambda^{2}}{2}\Expe{M_{\tau_{a} }F(a) }\ln\frac{1}{\e}\\
=&\frac{1}{\lambda}\para{\Psi(a,\lambda)-1}+\Expe{M_{\tau_{a} }\ln\para{M_{\tau_{a}}}F(a)}.
\end{eqalign}
The identity is
\begin{eqalign}
\Expe{M_{\tau_{a} }  }   =& 1+\lambda   \Expe{ M_{\tau_{a} } F(a)}.
\end{eqalign}
The derivative of the LHS is
\begin{eqalign}
\ptf{\Expe{M_{\tau_{a} }  } }{a}  =& \lambda\Expe{ M_{\tau_{a} } \ptf{U_{\e}(x)}{x}\mid_{x=\tau_{a}}\ptf{\tau_{a}}{a}}.
\end{eqalign}
The derivative of the RHS is
\begin{eqalign}
\ptf{\Expe{M_{\tau_{a} }  } }{a}  =& \lambda\Expe{ M_{\tau_{a} } \ptf{U_{\e}(x)}{x}\mid_{x=\tau_{a}}\ptf{\tau_{a}}{a} F(a)}\\
&+\lambda  \Expe{ M_{\tau_{a} } \drv{F(a)}{a} },
\end{eqalign}
where
\begin{eqalign}
 \drv{F(a)}{a}=&\drv{}{a}\para{\int_{0}^{a}    \ptf{R}{t}(\tau_{b},\tau_{a})\db }= \ptf{R}{t}(\tau_{b},\tau_{a})\mid_{b=a}+   \int_{0}^{a}    \ptf{^{2}R}{t_{1}\partial t_{2}}(\tau_{b},\tau_{a})\db\ptf{\tau_{a}}{a}.
\end{eqalign}

\end{enumerate}

\newpage\section{Moments of the maximum and minimum of modulus of GMC }\label{sec:maxminmodGMC}
In this section we study tail estimates and small ball estimates of the maximum/minimum of shifted GMC from \cite{binder2023inverse}.  One frequent theme is utilizing the 1d-correlation structure of GMC namely that neighboring evaluations $\eta[0,1],\eta[1,2],\eta[2,3],\eta[3,4]$ are correlated. But the pairs $\eta[0,1],\eta[2,3]$ and $\eta[1,2],\eta[3,4]$ are separately \iid. First we study the tail and moments of the maximum of the modulus of GMC.\\
On the face of it, in studying the $\maxl{0\leq T\leq L}\etamu{T,T+x}{\delta}$, we see that it could diverge as $\delta,x\to 0$ because we might be able to lower bound it by an increasing sequence of iid random variables such as $\etam{kx,x(k+1)}$ for $k\in [1,\floor{\frac{L}{x}}]$. We will see that at least for fixed $\delta>0$, we actually do have decay as $x\to 0$. This is in the spirit of chaining techniques where supremum over a continuum index set is dominated in terms of a maximum over a  finite index set.\\~\\
\noindent We will also need an extension for a different field: for $\lambda<1$, the field $U_{ \varepsilon}^{\delta, \lambda}$ with covariance
\begin{equation}
\Expe{U_{ \varepsilon}^{  \delta,\lambda }(x_{1} )U_{ \varepsilon}^{  \delta,\lambda }(x_{2} )  }=\left\{\begin{matrix}
\ln(\frac{\delta }{\varepsilon} )-\para{\frac{1}{\e}-\frac{1}{\delta}}\abs{x_{2}-x_{1}}+(1-\lambda)(1-\frac{\abs{x_{2}-x_{1}}}{\delta})&\tifc \abs{x_{2}-x_{1}}\leq \varepsilon\\ ~\\
 \ln(\frac{\delta}{\abs{x_{2}-x_{1}}})-1+\frac{\abs{x_{2}-x_{1}}}{\delta}+(1-\lambda)(1-\frac{\abs{x_{2}-x_{1}}}{\delta}) &\tifc \e\leq \abs{x_{2}-x_{1}}\leq \frac{\delta}{\lambda}\\~\\
  0&\tifc \frac{\delta}{\lambda}\leq \abs{x_{2}-x_{1}}
\end{matrix}\right.    .
\end{equation}
\begin{proposition}\label{prop:maxmoduluseta}
\pparagraph{Moments $p\in [1,\frac{2}{\gamma^{2}})$ }
For $L,\delta,x\geq 0$  and $\delta\leq 1$ we have
\begin{equation}\label{eq:maxmodulusetapone}
\Expe{\para{\supl{T\in[0,L] }\etamu{T,T+x}{\delta}}^{p}}\leq cx^{\alpha(p)}\para{\ceil{\frac{L}{x}+1}}^{\frac{p}{r_{p}}}\leq c(1+L+x)^{\frac{p}{r_{p}}} x^{\alpha(p)-\frac{p}{r_{p}}},    
\end{equation}
where $\alpha(p)=\zeta(p)$ when $x\leq 1$ and $\alpha(p)=p$ when $x\geq 1$, and the $r_{p}>0$ is an arbitrary number in $p<r_{p}<\frac{2}{\gamma^{2}}$. For simplification, we will also write $\frac{p}{r_{p}}=p(\frac{\gamma^{2}}{2}+\e_{p})$ for small enough $\e_{p}>0$. The same estimate follows for the measure $\eta^{\delta,\lambda}$ when $x\leq \delta$. 
\pparagraph{Moments $p\in (0,1)$}Here we have
\begin{eqalign}\label{eq:plessthaone}
\Expe{\para{\supl{T\in[0,L] }\etamu{T,T+x}{\delta}}^{p}}&\lessapprox \para{(1+L+x)^{\frac{1}{r_{1}}} x^{1-\frac{1}{r_{1}}}}^{p},    
\end{eqalign}
where as above $1<r_{1}<\frac{2}{\gamma^{2}}$ and let $c_{1}:=\frac{r_{1}-1}{r_{1}}=1-\beta-\e$ for arbitrarily small $\e>0$.
\end{proposition}
\begin{remark}
In \cref{eq:maxmodulusetapone}, we see that when  $\alpha(p)-\frac{p}{r_{p}}>0$, it decays to zero as $x\to 0$.   By taking $r_{p}\approx \frac{2}{\gamma^{2}}$, that means we require $\zeta(p)-\frac{p}{r_{p}}\approx p\frac{\gamma^{2}}{2}(\frac{2}{\gamma^{2}}-p)>0$. Also, one can check that this exponent is a bit better than that given in \cite[10.1 Theorem]{schilling2021brownian} for general stochastic processes.
\end{remark}
\noindent Next we study the negative moments for the minimum of the modulus of GMC. 
\begin{proposition}\label{prop:minmodeta}
We have for $p>0$
\begin{equation}
\Expe{\para{\infl{T\in[0,L] }\etamu{T,T+x}{\delta}}^{-p}}\lessapprox~x^{a_{\delta}(-p)}\para{\frac{L}{x}+2}^{\frac{p}{r}}2^{-\zeta(-r)\frac{p}{r}},
\end{equation}
where $a_{\delta}(-p):=\zeta(-p)$ when $x\leq \delta$ and $a_{\delta}(-p):=-p$ when $x\geq \delta$ and $r>0$ satisfies $\frac{p}{r}<1$ and so for simplicity we take arbitrarily small $\e_{p}:=\frac{p}{r}>0$.  The same follows for the measure $\eta^{\delta,\lambda}$ and $x\leq \delta$.
\end{proposition}
\begin{remark}
Here we note that as $r\to +\infty$, the constant  $2^{-\zeta(-r)\frac{p}{r}}$ diverges. So the smaller $\e_{p}:=\frac{p}{r}>0$, the larger the comparison constant.
\end{remark}

\newpage\section{Propeties of the covariance of truncated field}

\subsection{Regularity of the covariance}\label{sec:regcovariance}
The following are some of the hypotheses used in the development of Malliavin calculus for Gaussian processes \cite{decreusefond2008hitting,lei2012stochastic}.  The difference is
\begin{equation}
\Expe{\abs{U_{\varepsilon}^{r}(t)-U_{\varepsilon}^{r}(s)}^{2}}    =2\frac{\abs{t-s}}{\e}(1-\frac{\e}{r}),
\end{equation}
which is strictly positive for $t\neq s$. The covariance
\begin{equation}
R(\tau,t):= \left\{\begin{matrix}
\ln(\frac{r }{\varepsilon} )-\para{\frac{1}{\e}-\frac{1}{r}}\abs{\tau-t} &, \abs{\tau-t}\leq \varepsilon\\ 
 \ln(\frac{r}{\abs{\tau-t}}) +\frac{\abs{\tau-t}}{r}-1&, r>\abs{\tau-t}\geq \e
\end{matrix}\right.       
\end{equation}
is in fact an absolutely continuous function as a map $t\mapsto R(\tau,t)$ for each $\tau$: when $\abs{\tau-t}\leq \varepsilon$,  we have the absolutely continuous function $g(t)=\abs{\tau-t}$, and when $\abs{\tau-t}> \varepsilon$, we use that $\ln\frac{1}{x}$ is a differentiable function for $x>0$. We compute the partial derivative to be
\begin{equation}\label{eq:covariancederivative}
\ptf{R(\tau,t)}{t}= \left\{\begin{matrix}
-\para{\frac{1}{\e}-\frac{1}{r}}\frac{t-\tau}{\abs{t-\tau}} &, \abs{\tau-t}\leq \varepsilon\\ 
 -\frac{1}{\abs{t-\tau}}\frac{t-\tau}{\abs{t-\tau}} +\frac{1}{r}\frac{t-\tau}{\abs{t-\tau}}&,r>\abs{\tau-t}\geq \e
\end{matrix}\right.       .
\end{equation}
Therefore, for $t>\tau$ the derivative is negative $\ptf{R(\tau,t)}{t}<0$ and for $t<\tau$ it is positive $\ptf{R(\tau,t)}{t}>0$. So it is not continuous on the diagonal, which was one of the constraints in \cite{decreusefond2008hitting}. However, in the work  \cite{lei2012stochastic}, they manage to weaken to the following hypotheses that are satisfied here.
\begin{lemma}\label{lem:regularitycovariance}
For all $T>0$ the supremum of the integral of the partial derivative is finite for any $\alpha\geq 1$
\begin{equation}\label{finitint}
\sup_{s\in [0,T]}\int_{0}^{T}\abs{\ptf{R(s,t)}{t}}^{\alpha}  \dt<\infty  
\end{equation}
with a bound that diverges as $T\to +\infty$ or $\e\to 0$. In fact for any continuous function $f$ we have that
\begin{equation}\label{eq:continuousint}
s\mapsto F(s):=\int_{0}^{T}   f(t)\ptf{R(s,t)}{t}\dt  
\end{equation}
is continuous on $[0,\infty)$ as long as $\e>0$.
\end{lemma}
\begin{proofs}

\proofparagraph{Finite integral: proof of \cref{finitint}}
\pparagraph{Case $\alpha=1$}
Because for $\abs{s-t}\geq r$, we have zero covariance, we restrict the integral to the domains
\begin{equation}
[(s-r)\vee 0,(s-\e)\vee 0]    \cup [(s-\e)\vee 0,s]\cup [s,(s+\e)\wedge T]\cup [(s+\e)\wedge T,(s+r)\wedge T].
\end{equation}
In the domain $[(s-r)\vee 0,(s-\e)\vee 0]$, we have $t<s$ and $s-t>\e$ and so $\abs{\ptf{R(s,t)}{t}}=\frac{1}{s-t}-\frac{1}{r}$ and the integral will be
\begin{equation}
\int_{(s-r)\vee 0}^{(s-\e)\vee 0}\frac{1}{s-t}-\frac{1}{r}\dt=\ln(\frac{r\wedge s}{\e\wedge s})-\frac{1}{r} \para{s\wedge r-s\wedge \e} \end{equation}
Similarly, in the domain $ [(s+\e)\wedge T,(s+r)\wedge T]$, we have $\abs{\ptf{R(s,t)}{t}}=\abs{-\para{\frac{1}{t-s}-\frac{1}{r}}}=\frac{1}{t-s}-\frac{1}{r}$ and the integral will be
\begin{equation}
\ln(\frac{r\wedge (T-s)}{\e\wedge (T-s)})-\frac{1}{r} \para{(T-s)\wedge r-(T-s)\wedge \e}
\end{equation}
In the domain $[(s-\e)\vee 0,s]$, we have $\abs{\ptf{R(s,t)}{t}}=\para{\frac{1}{\e}-\frac{1}{r}}=:c_{\e,r}$ and similarly, in $[s,(s+\e)\wedge T]$ we again have $\abs{\ptf{R(s,t)}{t}}=\abs{-\para{\frac{1}{\e}-\frac{1}{r}}}=:c_{\e,r}$. Therefore, the total integral will be
\begin{equation}
\ln(\frac{r\wedge s}{\e\wedge s})-\frac{1}{r} \para{s\wedge r-s\wedge \e}+\ln(\frac{r\wedge (T-s)}{\e\wedge (T-s)})-\frac{1}{r} \para{(T-s)\wedge r-(T-s)\wedge \e}+c_{\e,r}\para{(s+\e)\wedge T-(s-\e)\vee 0}.
\end{equation}
So we see from here that as $\e\to 0$, this integral diverges. The log-terms are the only source of potential singularity. When $s$ is close to zero i.e. $r>s>\e$ or $\e\geq s$, we get $\ln(\frac{s}{\e})$ and  $\ln(\frac{s}{s})=0$ respectively. When $s$ is close to $T$ i.e. $r>T-s>\e$ or $\e\geq T-s$, we similarly get
$\ln(\frac{T-s}{\e})$ and  $\ln(\frac{T-s}{T-s})=0$ respectively. Therefore, we indeed have a finite supremum for each $T>0$. 
\pparagraph{Case $\alpha>1$} 
Here instead of logarithms we get singular terms of the form $\frac{1}{x^{\alpha-1}}$. In particular following the same integration steps on splitting domains we get singular terms of the following form:
\begin{equation}
\frac{1}{(r\wedge s)^{\alpha-1}}-\frac{1}{(\e\wedge s)^{\alpha-1}}    \tand\frac{1}{(r\wedge (T-s))^{\alpha-1}}-\frac{1}{(\e\wedge  (T-s))^{\alpha-1}}.
\end{equation}
When $s$ is close to zero i.e. $r>s>\e$ or $\e\geq s$, we get $\frac{1}{r^{\alpha-1}}-\frac{1}{\e^{\alpha-1}}$ and  $\frac{1}{s^{\alpha-1}}-\frac{1}{s^{\alpha-1}}=0$ respectively. For $s$ close to $T$, we conversely get  $\frac{1}{r^{\alpha-1}}-\frac{1}{\e^{\alpha-1}}$ and  $\frac{1}{(T-s)^{\alpha-1}}-\frac{1}{(T-s)^{\alpha-1}}=0$. We always get a singular power in $\e>0$. In summary, we again have a finite supremum for each $T>0$ and $\e>0$. 
\proofparagraph{The continuous weighted derivative: proof of \cref{eq:continuousint}}
We split over the same domains. We end up with the following total integral
\begin{eqalign}
 &\int_{(s-r)\vee 0}^{(s-\e)\vee 0}\frac{f(t)}{s-t}\dt+ (-\frac{1}{r})   \int_{(s-r)\vee 0}^{(s-\e)\vee 0}f(t)\dt+ \int_{(s+\e)\wedge T}^{(s+r)\wedge T}\frac{f(t)}{t-s}\dt+ (-\frac{1}{r})   \int_{(s+\e)\wedge T}^{(s+r)\wedge T}f(t)\dt\\
 &+c_{\e,r} \int_{(s-\e)\vee 0}^{(s+\e)\wedge T}f(t)\dt. 
\end{eqalign}
The integrals containing only the continuous function $f(t)$ are differentiable in $s$ due to the fundamental theorem of calculus. In particular, the function $g(t)=\frac{1}{s-t}$ is continuously differentiable in the above domains because they don't contain an $\e$-neighbourhood of the singularity $t=s$. Therefore, the integrals with integrands $\frac{f(t)}{s-t}$ are differentiable due to Leibniz-rule. 
\pparagraph{Case of $\e\to 0$ and large $T$}
Here we get
\begin{eqalign}
 &\int_{(s-r)\vee 0}^{s}\frac{f(t)}{s-t}\dt+ (-\frac{1}{r})   \int_{(s-r)\vee 0}^{s}f(t)\dt+ \int_{s}^{s+r}\frac{f(t)}{t-s}\dt+ (-\frac{1}{r})   \int_{s}^{s+r}f(t)\dt\\
 &+\liz{\e}\para{\frac{1}{\e}-\frac{1}{r}} \int_{(s-\e)\vee 0}^{(s+\e)\wedge T}f(t)\dt. 
\end{eqalign}

\end{proofs}



\printbibliography[title={Whole bibliography}]

\end{document}